\documentclass{article}
\begin{document}
\font\germ=eufm10
\def\ssl{\hbox{\germ sl}}
\def\slh{\widehat{\ssl_2}}
\makeatletter
\def\aaa{@}
\centerline{}
\vskip1cm
\centerline{\Large\bf Polytopes for  Crystallized Demazure Modules and }
\centerline{\Large\bf Extremal Vectors}
\vskip5pt
\centerline{Toshiki NAKASHIMA}
\vskip4pt
\centerline{Department of Mathematics,}
\centerline{Sophia University, Tokyo 102-8554, JAPAN}
\centerline{e-mail:\,\,toshiki@mm.sophia.ac.jp}
\makeatother

\renewcommand{\labelenumi}{$($\roman{enumi}$)$}
\renewcommand{\labelenumii}{$(${\rm \alph{enumii}}$)$}
\font\germ=eufm10

\def\al{\alpha}
\def\beneme{\begin{enumerate}}
\def\beq{\begin{equation}}
\def\beqn{\begin{eqnarray}}
\def\beqnn{\begin{eqnarray*}}
\def\bigsl{{\hbox{\fontD \char'54}}}
\def\cd{\cdots}
\def\ddd{\hbox{\germ D}}
\def\del{\delta}
\def\Del{\Delta}
\def\ei{e_i}
\def\eit{\tilde{e}_i}
\def\eneme{\end{enumerate}}
\def\ep{\epsilon}
\def\eeq{\end{equation}}
\def\eeqn{\end{eqnarray}}
\def\eeqnn{\end{eqnarray*}}
\def\fit{\tilde{f}_i}
\def\ft{\tilde{f}}
\def\ge{\hbox{\germ g}}
\def\gl{\hbox{\germ gl}}
\def\hom{{\hbox{Hom}}}
\def\ify{\infty}
\def\io{\iota}
\def\kp{k^{(+)}}
\def\km{k^{(-)}}
\def\llra{\relbar\joinrel\relbar\joinrel\relbar\joinrel\rightarrow}
\def\lan{\langle}
\def\lar{\longrightarrow}
\def\lm{\lambda}
\def\Lm{\Lambda}
\def\mapright#1{\smash{\mathop{\longrightarrow}\limits^{#1}}}
\def\nd{\noindent}
\def\nn{\nonumber}
\def\nnn{\hbox{\germ n}}
\def\oint{{\cal O}_{\rm int}(\ge)}
\def\ot{\otimes}
\def\op{\oplus}
\def\opi{\ovl\pi_{\lm}}
\def\ovl{\overline}
\def\plm{\Psi^{(\lm)}_{\io}}
\def\qq{\qquad}
\def\q{\quad}
\def\qed{\hfill\framebox[2mm]{}}
\def\QQ{\hbox{\bf Q}}
\def\qi{q_i}
\def\qii{q_i^{-1}}
\def\ran{\rangle}
\def\rlm{r_{\lm}}
\def\ssl{\hbox{\germ sl}}
\def\slh{\widehat{\ssl_2}}
\def\ti{t_i}
\def\tii{t_i^{-1}}
\def\til{\tilde}
\def\tt{{\hbox{\germ{t}}}}
\def\ttt{\hbox{\germ t}}
\def\uq{U_q(\ge)}
\def\uqm{U^-_q(\ge)}
\def\uqp{U^+_q(\ge)}
\def\uqmq{{U^-_q(\ge)}_{\bf Q}}
\def\uqpm{U^{\pm}_q(\ge)}
\def\uqq{U_{\bf Q}^-(\ge)}
\def\uqz{U^-_{\bf Z}(\ge)}
\def\util{\widetilde\uq}
\def\vep{\varepsilon}
\def\vp{\varphi}
\def\vpi{\varphi^{-1}}
\def\xii{\xi^{(i)}}
\def\Xiioi{\Xi_{\io}^{(i)}}
\def\wtil{\widetilde}
\def\what{\widehat}
\def\wpi{\widehat\pi_{\lm}}
\def\ZZ{\hbox{\bf Z}}

\renewcommand{\thesection}{\arabic{section}}
\section{Introduction}
\setcounter{equation}{0}
\renewcommand{\theequation}{\thesection.\arabic{equation}}

Demazure's character formula for arbitrary Kac-Moody Lie algebra
was given by S.Kumar and O.Mathieu independently (\cite{KS},\cite{M})
by using geometric methods.
In 1995, P.Littelmann gave some conjecture (partially solved by
himself) about the 
relation between Demazure's character formula and crystal bases \cite{L},
which was solved affirmatively by M.Kashiwara \cite{K3}.
Then it gave purely algebraic proof for Demazure's character
formula for symmetrizable Kac-Moody Lie algebras.
Here let us see those formulations.
Let $\ge$ be a symmetrizable Kac-Moody Lie algebra
(in the context of ``crystal base'', we need ``symmetrizable''),
and $\nnn^+$ be the nilpotent subalgebra of $\ge$.
Furthermore, let $\ZZ[P]$ be the group algebra of the 
weight lattice $P$ and
$W$ be the Weyl group associated with $\ge$.
Then Demazure operator 
$D_w:\ZZ[P]\longrightarrow \ZZ[P]$ ($w\in W$)
is given as follows: for $i\in I$ (index set) we set
$D_i(e^{\lm}):=e^{\lm}(1-e^{-{(1+\lan h_i,\lm\ran)\al_i}})/1-e^{-\al_i}$
and for $w=s_{i_l}\cd s_{i_1}$ set $D_w:=D_{i_l}\cd D_{i_1}$,
which is well-defined.
Let $V(\lm)$ be the irreducible highest weight module with the 
highest weight $\lm$ and $u_{w\lm}$ be the 
extrmal vector with the weight $w\lm$ $(w\in W)$.
Then, Demazure's character formula is described as follows:
\begin{equation}
\label{DCF}
ch(U(\nnn^+)u_{w\lm})=D_w(e^{\lm}).
\end{equation}
In \cite{L}, Littelmann gave the following conjecture:
Let $V(\lm)$ be the irreducible $\uq$-highest weight module with the 
highest weight $\lm$ and $(L(\lm),B(\lm))$ be its crystal base.
Then there exists a subset $B_w(\lm)\subset B(\lm)$ such that 
\begin{eqnarray}
&&\uqp u_{w\lm}\cap L(\lm)/{\uqp u_{w\lm}\cap qL(\lm)}
=\bigoplus _{b\in B_w(\lm)}\QQ b,
\label{L1}\\
&& \sum_{b\in B_w(\lm)}b={\ddd}_{i_l}\cd{\ddd}_{i_1}u_{\lm},
\label{L2}
\end{eqnarray}
where $u_{\lm}$ is the highest weight vector
with the weight $\lm$ and 
${\ddd}_i$ is the additive operator on $\ZZ^{\oplus B(\lm)}$
given by:
$$
{\ddd}_i\, b:=\left\{
\begin{array}{ll}
\sum_{0\leq k\leq \lan h_i,wt(b)\ran}\fit^k b &
{\mbox if }\lan h_i,wt(b)\ran\geq0,\\
-\sum_{1 \leq k< -\lan h_i,wt(b)\ran}\eit^k b &
{\mbox if }\lan h_i,wt(b)\ran<0.
\end{array}
\right.
$$
We call the left-hand side of (\ref{L1}) 
{\it crystallized Demazure module} of $V(\lm)$ 
assciated with $w\in W$.
Here we know that Littelmann's conjecture implies 
Demazure's character formula by the following way:
Define the operator $ewt:\ZZ^{\oplus B(\lm)}\longrightarrow \ZZ[P]$
by $ewt(b):=e^{wt(b)}$ for $b\in B(\lm)$ and 
$ewt(b_1+b_2)=ewt(b_1)+ewt(b_2)$. Now, 
we have $ewt({\ddd}_ib)=D_i(ewt(b)).$
Thus, by (\ref{L1}) and (\ref{L2}) we have 
$$
\begin{array}{l}
ch(\uqp u_{w\lm})=ewt(\sum_{b\in B_w(\lm)}b)
=ewt({\ddd}_{i_l}\cd{\ddd}_{i_1}u_{\lm})\\
=D_{i_l}\cd D_{i_1} ewt(u_{\lm})=
D_{i_l}\cd D_{i_1}(e^{\lm})=D_w(e^{\lm}).
\end{array}
$$
In \cite{K3}, Kashiwara shown the existence of $B_w(\lm)$ 
for arbitrary symmetrizable Kac-Moody cases and characterized it as follows:
\newtheorem{thm1}{Theorem}[section]
\begin{thm1}[\cite{K3}]
\label{kas}
\begin{enumerate}
\item
$\eit B_w(\lm)\subset B_w(\lm)\sqcup\{0\}.$
\item
If $s_i w<w$ (Bruhat order), then 
$B_w(\lm)=\{\fit^k b;k\geq 0,\,\,b\in
B_{s_iw}(\lm),\,\,\eit b=0\}\setminus\{0\}.$
\item
For any $i$-string $S$, $S\cap B_w(\lm)$ is either 
empty or $S$ or $\{$the highest weight vector of $S$ $\}$.
\end{enumerate}
\end{thm1}

In \cite{NZ},\cite{N2}, we developed the polyhedral realization
of crystal bases.
We shall explain the relations between crystal bases 
of Demazure modules and the polyhedral realizations briefly.
Let $\io=\cd i_k,\cd,i_2,i_1$ be an infinite sequence from the index
set $I$ satisfying some condition and $\lm$ be a dominant 
integral weight. Then 
there exists the embedding $\Psi^{(\lm)}_\io:B(\lm)\hookrightarrow
\ZZ^{\ify}_{\io}[\lm](\cong \ZZ^{\ify})$. 
The exact image of $\Psi^{(\lm)}_{\io}$ is described 
(under some assumption)
as a subset in $\ZZ^{\ify}$ 
given by some system of linear inequalities,
which is called {\it polyhedral realization}.
Let $w=s_{i_l}\cd s_{i_1}$ (reduced expression) be an element in $W$
and take a sequence $\io=(j_k)_{k\geq 1}$ which satisfies 
$i_k=j_k$ ($1\leq k\leq l$).
Then in this paper,
the subset $\Psi^{(\lm)}_{\io}(B_w(\lm))$ is given as a set of 
lattice points of some convex polytope in $\ZZ^{\ify}$,
where ``polytope'' means a bounded polyhedron.
Furthermore, we succeed in giving explicit form of 
extremal vector $\Psi^{(\lm)}_{\io}(u_{w\lm})$
which is contained in $\Psi^{(\lm)}_{\io}(B_w(\lm))$
as the  unique solution of some system of linear equations.

The organization of this paper is as follows:
In Sect.2 we review the polyhedral
realizations of crystals.
We shall describe the polytopes for 
$B_w(\lm)$ in Sect. 3 and 
the extremal vectors in Sect.4.

\section{Polyhedral realizations of crystals}
\setcounter{equation}{0}
\renewcommand{\thesection}{\arabic{section}}
\renewcommand{\theequation}{\thesection.\arabic{equation}}

\subsection{Notations}

We list the notations used in this paper.
Most of them are same as those in \cite{N2}.

Let $\ge$ be
a  symmetrizable Kac-Moody algebra over {\bf Q}
with a Cartan subalgebra
$\ttt$, a weight lattice $P \subset \ttt^*$, the set of simple roots
$\{\al_i: i\in I\} \subset \ttt^*$,
and the set of coroots $\{h_i: i\in I\} \subset \ttt$,
where $I$ is a finite index set.
Let $\lan h,\lm\ran$ be the pairing between $\ttt$ and $\ttt^*$,
and $(\al, \beta)$ be an inner product on
$\ttt^*$ such that $(\al_i,\al_i)\in 2{\bf Z}_{\geq 0}$ and
$\lan h_i,\lm\ran={{2(\al_i,\lm)}\over{(\al_i,\al_i)}}$
for $\lm\in\ttt^*$.
Let $P^*=\{h\in \ttt: \lan h,P\ran\subset\ZZ\}$ and
$P_+:=\{\lm\in P:\lan h_i,\lm\ran\in\ZZ_{\geq 0}\}$.
We call an element in $P_+$ a {\it dominant integral weight}.
Here we define a partial order on $P$ by:
For $\lm,\mu\in P$, $\lm\succ \mu$ $\Leftrightarrow$
$\lm-\mu\in \oplus_{i\in I}\QQ_{\geq0}\al_i$.
The quantum algebra $\uq$
is an associative
$\QQ(q)$-algebra generated by the $e_i$, $f_i \,\, (i\in I)$,
and $q^h \,\, (h\in P^*)$
satisfying the usual relations.
The algebra $\uqm$ is the subalgebra of $\uq$ generated 
by the $f_i$ $(i\in I)$.

For the irreducible highest weight module of $\uq$
with the highest weight $\lm\in P_+$, we denote  $V(\lm)$
and its {\it crystal base} we denote $(L(\lm),B(\lm))$.
Similarly, for the crystal base of the algebra $\uqm$ we denote 
$(L(\ify),B(\ify))$ (see \cite{K0},\cite{K1},\cite{K6}).
Let $\pi_{\lm}:\uqm\longrightarrow V(\lm)\cong \uqm/{\sum_i\uqm
\fit^{1+\lan h_i,\lm\ran}}$ be the canonical projection and 
$\widehat \pi_{\lm}:L(\ify)/qL(\ify)\longrightarrow L(\lm)/qL(\lm)$
be the induced map from $\pi_{\lm}$. Here note that 
$\widehat \pi_{\lm}(B(\ify))=B(\lm)\sqcup\{0\}$.

By the terminology {\it crystal } we mean some combinatorial object 
obtained by abstracting the properties of crystal bases.
Indeed, crystal constitutes a set $B$ and the maps
$wt:B\longrightarrow P$, $\vep_i,\vp_i:B\longrightarrow \ZZ\sqcup\{-\ify\}$
and $\eit,\fit:B\sqcup\{0\}\longrightarrow B\sqcup\{0\}$
($i\in I$) with several axioms (see \cite{K3},\cite{NZ},\cite{N2}).
In fact, $B(\ify)$ and $B(\lm)$ are the typical examples 
of crystals.

It is well-known that $\uq$ has a Hopf algebra structure.
Then the tensor product of $\uq$-modules has
a $\uq$-module structure.
The crystal bases have very nice properties for 
tensor operations. Indeed, if $(L_i,B_i)$ is a crystal base of 
$\uq$-module $M_i$ ($i=1,2$), $(L_1\ot_A L_2, B_1\ot B_2)$
is a crystal base of $M_1\ot_{\QQ(q)} M_2$ (\cite{K1}).
Consequently, we can consider the tensor product
of crystals and then they constitute a tensor category.

\subsection{Polyhedral Realization of $B(\ify)$}
\label{poly-uqm}
In this subsection, we recall the results in \cite{NZ}.

%
Consider the additive group
\begin{equation}
\ZZ^{\ify}
:=\{(\cd,x_k,\cd,x_2,x_1): x_k\in\ZZ
\,\,{\rm and}\,\,x_k=0\,\,{\rm for}\,\,k\gg 0\};
\label{uni-cone}
\end{equation}
we will denote by $\ZZ^{\ify}_{\geq 0} \subset \ZZ^{\ify}$
the subsemigroup of nonnegative sequences.
To the rest of this section, we fix an infinite sequence of indices
$\io=\cd,i_k,\cd,i_2,i_1$ from $I$ such that
\begin{equation}
{\hbox{
$i_k\ne i_{k+1}$ and $\sharp\{k: i_k=i\}=\ify$ for any $i\in I$.}}
\label{seq-con}
\end{equation}

We can associate to $\io$ a crystal structure
on $\ZZ^{\ify}$ and denote it by $\ZZ^{\ify}_{\io}$ 
(\cite[2.4]{NZ}).

\newtheorem{pro2}{Proposition}[section]
\begin{pro2}[\cite{K3}, See also \cite{NZ}]
\label{emb}
There is a unique embedding of crystals
$($called Kashiwara embedding$)$
\begin{equation}
\Psi_{\io}:B(\ify)\hookrightarrow \ZZ^{\ify}_{\geq 0}
\subset \ZZ^{\ify}_{\io},
\label{psi}
\end{equation}
such that
$\Psi_{\io} (u_{\ify}) = (\cd,0,\cd,0,0)$.
\end{pro2}

Consider the infinite dimensional vector space
$$
\QQ^{\ify}:=\{{x}=
(\cd,x_k,\cd,x_2,x_1): x_k \in \QQ\,\,{\rm and }\,\,
x_k = 0\,\,{\rm for}\,\, k \gg 0\},
$$
and its dual space $(\QQ^{\ify})^*:={\rm Hom}(\QQ^{\ify},\QQ)$.
We will write a linear form $\vp \in (\QQ^{\ify})^*$ as
$\vp({x})=\sum_{k \geq 1} \vp_k x_k$ ($\vp_j\in \QQ$).

For the fixed infinite sequence
$\io=(i_k)$ we set $\kp:={\rm min}\{l:l>k\,\,{\rm and }\,\,i_k=i_l\}$ and
$\km:={\rm max}\{l:l<k\,\,{\rm and }\,\,i_k=i_l\}$ if it exists,
or $\km=0$  otherwise.
We set for $x\in \QQ^{\ify}$, $\beta_0(x)=0$ and
\begin{equation}
\beta_k(x):=x_k+\sum_{k<j<\kp}\lan h_{i_k},\al_{i_j}\ran x_j+x_{\kp}
\qq(k\geq1).
\label{betak}
\end{equation}
We define a piecewise-linear operator $S_k=S_{k,\io}$ on $(\QQ^{\ify})^*$ by
$$
S_k(\vp):=
\left\{
\begin{array}{ll}
\vp-\vp_k\beta_k & {\mbox{ if }}\vp_k>0,\\
 \vp-\vp_k\beta_{\km} & {\mbox{ if }}\vp_k\leq 0.
\end{array}
\right.
$$
Here we set
\begin{eqnarray}
\Xi_{\io} &:=  &\{S_{j_l}\cd S_{j_2}S_{j_1}x_{j_0}\,|\,
l\geq0,j_0,j_1,\cd,j_l\geq1\},
\label{Xi_io}\\
\Sigma_{\io} & := &
\{x\in \ZZ^{\ify}\subset \QQ^{\ify}\,|\,\vp(x)\geq0\,\,{\rm for}\,\,
{\rm any}\,\,\vp\in \Xi_{\io}\}.
\end{eqnarray}
We impose on $\io$ the following positivity assumption:
\begin{equation}
{\hbox{if $\km=0$ then $\vp_k\geq0$ for any 
$\vp(x)=\sum_k\vp_kx_k\in \Xi_{\io}$}}.
\label{posi}
\end{equation}
\newtheorem{thm2}{Theorem}[section]
\begin{thm2}[\cite{NZ}]
Let $\io$ be a sequence of indices satisfying $(\ref{seq-con})$ 
and (\ref{posi}). Then we have 
${\rm Im}(\Psi_{\io})(\cong B(\ify))=\Sigma_{\io}$.
\end{thm2}

\subsection{Polyhedral Realization of $B(\lm)$}

In this subsection, 
we review the result in \cite{N2}. 
In the rest of this section,
$\lm$ is supposed to be a dominant integral weight.
Let $R_{\lm}:=\{r_{\lm}\}$ be the crystal defined in \cite{N2}.
Consider the crystal $B(\ify)\ot R_{\lm}$ and  define the map
\begin{equation}
\Phi_{\lm}:(B(\ify)\ot R_{\lm})\sqcup\{0\}\longrightarrow B(\lm)\sqcup\{0\},
\label{philm}
\end{equation}
by $\Phi_{\lm}(0)=0$ and $\Phi_{\lm}(b\ot r_{\lm})=\wpi(b)$ for $b\in B(\ify)$.
We set
$$
\wtil B(\lm):=
\{b\ot r_{\lm}\in B(\ify)\ot R_{\lm}\,|\,\Phi_{\lm}(b\ot r_{\lm})\ne 0\}.
$$

\begin{thm2}[\cite{N2}]
\label{ify-lm}
\begin{enumerate}
\item
The map $\Phi_{\lm}$ becomes a surjective strict morphism of crystals
$B(\ify)\ot R_{\lm}\longrightarrow B(\lm)$.
\item
$\wtil B(\lm)$ is a subcrystal of $B(\ify)\ot R_{\lm}$, 
and $\Phi_{\lm}$ induces the
isomorphism of crystals $\wtil B(\lm)\mapright{\sim} B(\lm)$.
\end{enumerate}
\end{thm2}

Let us denote $\ZZ^{\ify}_{\io}\ot R_{\lm}$ by 
$\ZZ^{\ify}_{\io}[\lm]$. Here note that
since the crystal $R_{\lm}$ has  only one element,
as a set we can identify $\ZZ^{\ify}_{\io}[\lm]$ with
$\ZZ^{\ify}_{\io}$ but their crystal structures are different.
As for the explicit crystal structure of $\ZZ^{\ify}_{\io}[\lm]$,
see \ref{ZZ-io-lm} below.
By Theorem \ref{ify-lm}, we have the strict embedding of crystals
$\Omega_{\lm}:B(\lm)(\cong \wtil B(\lm))\hookrightarrow B(\ify)\ot R_{\lm}.$
Combining $\Omega_{\lm}$ and the
Kashiwara embedding $\Psi_{\io}$,
we obtain the following:
\begin{thm2}[\cite{N2}]
\label{embedding}
There exists the unique  strict embedding of crystals
\begin{equation}
\Psi_{\io}^{(\lm)}:B(\lm)\stackrel{\Omega_{\lm}}{\hookrightarrow}
B(\ify)\ot R_{\lm}
\stackrel{\Psi_{\io}\ot {\rm id}}{\hookrightarrow}
\ZZ^{\ify}_{\io}\ot R_{\lm}=:\ZZ^{\ify}_{\io}[\lm],
\label{Psi-lm}
\end{equation}
such that $\Psi^{(\lm)}_{\io}(u_{\lm})=(\cd,0,0,0)\ot r_{\lm}$.
\end{thm2}

\vskip5pt


We fix a sequence of indices $\io$
satisfying (\ref{seq-con}) and take a dominant integral weight 
$\lm\in P_+$.
For $k\geq1$ let $k^{(\pm)}$ be  the ones in \ref{poly-uqm}.
Let $\beta_k^{(\pm)}(x)$ be linear functions given by
\begin{eqnarray}
&& \q
\beta_k^{(+)} (x)  =  \sigma_k (x) - \sigma_{\kp} (x)
= x_k+\sum_{k<j<\kp}\lan h_{i_k},\al_{i_j}\ran x_j+x_{\kp},
\label{beta}\\
&&  \beta_k^{(-)} (x) 
\label{beta--}
\\
&& =
\left\{
\begin{array}{ll}
\sigma_{\km} (x) - \sigma_k (x)
=x_{\km}+\sum_{\km<j<k}\lan h_{i_k},\al_{i_j}\ran x_j+x_k &
 \hspace{-10pt} {\mbox{ if }}\km>0,\\
\sigma_0^{(i_k)} (x) - \sigma_k (x)
=-\lan h_{i_k},\lm\ran+\sum_{1\leq j<k}\lan h_{i_k},\al_{i_j}\ran x_j+x_k
&
 \hspace{-10pt}{\mbox{ if }}\km=0,
\end{array}
\right.\nn
\end{eqnarray}
(As for the functions $\sigma_k$ and $\sigma^{(i)}_0$
see (\ref{sigma}) and (\ref{sigma0}) below.).
Here note that
$\beta_k^{(+)}=\beta_k$ and 
$\beta_k^{(-)}=\beta_{\km}  {\hbox{ \,\,if\,\, $\km>0$}}$.

Using this notation, for every $k \geq 1$, we define 
an operator
$\what S_k = \what S_{k,\io}$ for a linear function 
$\vp(x)=c+\sum_{k\geq 1}\vp_kx_k$
$(c,\vp_k\in\QQ)$ on $\QQ^{\ify}$ by:

$$
\what S_k\,(\vp) :=\left\{
\begin{array}{ll}
\vp - \vp_k \beta_k^{(+)} & {\mbox{ if }}\vp_k > 0,\\
\vp - \vp_k \beta_k^{(-)} & {\mbox{ if }}\vp_k \leq 0.
\end{array}
\right.
$$

For the fixed sequence $\io=(i_k)$, 
in case $\km=0$ for $k\geq1$, there exists unique $i\in I$ such that $i_k=i$.
We denote such $k$ by $\io^{(i)}$, namely, $\io^{(i)}$ is the first
number $k$ such that $i_k=i$.
Here for $\lm\in P_+$ and $i\in I$ we set
\begin{equation}
\lm^{(i)}(x):=
-\beta^{(-)}_{\io^{(i)}}(x)=\lan h_i,\lm\ran-\sum_{1\leq j<\io^{(i)}}
\lan h_i,\al_{i_j}\ran x_j-x_{\io^{(i)}}.
\label{lmi}
\end{equation}

For $\io$ and a dominant integral weight $\lm$,
let $\Xi_{\io}[\lm]$ be the set of all linear functions
generatd by $\what S_k=\what S_{k,\io}$ 
from the functions $x_j$ ($j\geq1$)
and $\lm^{(i)}$ ($i\in I$), namely,
\begin{equation}
\begin{array}{ll}
\Xi_{\io}[\lm]&:=\{\what S_{j_l}\cd\what S_{j_1}x_{j_0}\,
:\,l\geq0,\,j_0,\cd,j_l\geq1\}
\\
&\cup\{\what S_{j_k}\cd \what S_{j_1}\lm^{(i)}(x)\,
:\,k\geq0,\,i\in I,\,j_1,\cd,j_k\geq1\}.
\end{array}
\label{Xi}
\end{equation}
Now we set
\begin{equation}
\Sigma_{\io}[\lm]
:=\{x\in \ZZ^{\ify}_{\io}[\lm](\subset \QQ^{\ify})\,:\,
\vp(x)\geq 0\,\,{\rm for \,\,any }\,\,\vp\in \Xi_{\io}[\lm]\}.
\label{Sigma}
\end{equation}

For a sequence $\io$ and a domiant integral weight $\lm$, a pair
$(\io,\lm)$ is called {\it ample}
if $\Sigma_{\io}[\lm]\ni\vec 0=(\cd,0,0)$.

\begin{thm2}[\cite{N2}]
\label{main}
Suppose that $(\io,\lm)$ is ample.
Then we have
${\rm Im}(\plm)(\cong B(\lm))=\Sigma_{\io}[\lm]$.
\end{thm2}

\section{Crystallized Demazure modules}
\setcounter{equation}{0}
\renewcommand{\thesection}{\arabic{section}}
\renewcommand{\theequation}{\thesection.\arabic{equation}}

\subsection{Structure of $\ZZ^{\ify}_{\io}[\lm]$}
\label{ZZ-io-lm}

We shall review an explicit crystal structure of
$\ZZ^{\ify}[\lm]$ in \cite{N2}.
Fix a sequence of indices $\io:=(i_k)_{k\geq 1}$ satisfying the condition
(\ref{seq-con}) and a weight $\lm\in P$.
(Here we do not necessarily assume that 
$\lm$ is dominant.)
As we stated before, 
we can identify $\ZZ^{\ify}$ with $\ZZ^{\ify}[\lm]$
as a set. Thus $\ZZ^{\ify}[\lm]$ can be regarded 
as a subset of $\QQ^{\ify}$, and then 
we denote an element in $\ZZ^{\ify}[\lm]$
by $ x=(\cd,x_k,\cd,x_2,x_1)$.
For $ x=(\cd,x_k,\cd,x_2,x_1)\in \QQ^{\ify}$
we define the linear functions
\begin{eqnarray}
\sigma_k(x)&:= &x_k+\sum_{j>k}\lan h_{i_k},\al_{i_j}\ran x_j,
\q(k\geq1)
\label{sigma}\\
\sigma_0^{(i)}(x)
&:= &-\lan h_i,\lm\ran+\sum_{j\geq1}\lan h_i,\al_{i_j}\ran x_j,
\q(i\in I)
\label{sigma0}
\end{eqnarray}
Here note that
since $x_j=0$ for $j\gg0$ on $\QQ^{\ify}$,
the functions $\sigma_k$ and $\sigma^{(i)}_0$ are
well-defined.
Let $\sigma^{(i)} (x)
 := {\rm max}_{k: i_k = i}\sigma_k (x)$, and
$M^{(i)} :=
\{k: i_k = i, \sigma_k (x) = \sigma^{(i)}(x)\}.
$
Note that
$\sigma^{(i)} (x)\geq 0$, and that
$M^{(i)} = M^{(i)} (x)$ is a finite set
if and only if $\sigma^{(i)} (x) > 0$.
Now we define the maps
$\eit: \ZZ^{\ify}[\lm] \sqcup\{0\}\lar \ZZ^{\ify}[\lm] \sqcup\{0\}$
and
$\fit: \ZZ^{\ify}[\lm] \sqcup\{0\}\lar \ZZ^{\ify}[\lm] \sqcup\{0\}$ 
by setting $\eit(0)=\fit(0)=0$, and 
\begin{equation}
(\fit(x))_k  = x_k + \delta_{k,{\rm min}\,M^{(i)}}
\,\,{\rm if }\,\,\sigma^{(i)}(x)>\sigma^{(i)}_0(x);
\,\,{\rm otherwise}\,\,\fit(x)=0,
\label{action-f}
\end{equation}
\begin{equation}
(\eit(x))_k  = x_k - \delta_{k,{\rm max}\,M^{(i)}} \,\, {\rm if}\,\,
\sigma^{(i)} (x) > 0\,\,
{\rm and}\,\,\sigma^{(i)}(x)\geq\sigma^{(i)}_0(x) ; \,\,
 {\rm otherwise} \,\, \eit(x)=0,
\label{action-e}
\end{equation}
where $\del_{i,j}$ is the Kronecker's delta.
We also define the functions
$wt$, $\vep_i$ and $\vp_i$ on $\ZZ^{\ify}[\lm]$ by
\begin{eqnarray}
&& wt(x) :=\lm -\sum_{j=1}^{\ify} x_j \al_{i_j},
\label{wt-vep-vp-1}\\
&& \vep_i (x) := {\rm max}(\sigma^{(i)} (x),\sigma^{(i)}_0(x))
\label{wt-vep-vp}\\
&& \vp_i (x) := \lan h_i, wt(x) \ran + \vep_i(x).
\label{wt-vep-vp-3}
\end{eqnarray}
Note that 
by (\ref{wt-vep-vp-1}) we have
\begin{equation}
\lan h_i,wt(x)\ran = -\sigma^{(i)}_0(x).
\label{**}
\end{equation}
\subsection{Polytopes for $B_w(\lm)$}

In this section, we describe the explicit form of 
the polytopes corresponding to the crystals of Demazure module 
$B_w(\lm)$ $(\lm\in P_+)$ as in the introduction.

By the characterization of $B_w(\lm)$ given in Theorem \ref{kas} (ii),
we can construct it inductively according to some
reduced expression of $w$.
Indeed, we have $B_1(\lm)=\{u_{\lm}\}$ (where 1 is the 
identity of $W$) and then 
we obtain $B_{s_i}(\lm)=\{\fit^ku_{\lm}; k\geq0\}\setminus\{0\}$.
Since $\vec 0:=(\cd,0,0)$ corresponds to the highest weight vector,
by (\ref{action-f}) the image by $\Psi^{(\lm)}_{\io}$ is given by;
$$
\Sigma_{s_i}[\lm]:=\{(\cd,0,0,k);0 \leq k\leq \lan h_i,\lm\ran\},
$$
where $i_1=i$ for $\io=(i_k)_{k\geq1}$.

For $w\in W$, let us fix one reduced expression
$w=s_{i_L}s_{i_{L-1}}\cd s_{i_2}s_{i_1}$ and 
let $\io:=(j_k)_{k\geq1}$ be the infinite sequence of 
indices such that $i_k=j_k$ for $1\leq k\leq L$.
Here we do not necessarily assume that $(\io,\lm)$ is ample.
In this setting, we have

\newtheorem{pro3}{Proposition}[section]
\begin{pro3}
\label{polytope}
Set 
\begin{equation}
\label{sigma-w}
\Sigma_w[\lm]:=\{(\cd,x_k,x_{k-1},\cd,x_2,x_1)\in
{\rm Im}(\Psi^{(\lm)}_{\io})\,|\,x_k=0\,\,{\mbox{ for }}k>L\}.
\end{equation}
Then we have $\Psi_{\io}^{(\lm)}(B_w(\lm))=\Sigma_w[\lm]$.
\end{pro3}

{\sl Proof.\,\,}
We shall show by induction on the length of $w$.
If the length of $w$ is equal to 0, then $w=1$.
So we have $B_1(\lm)=\{u_{\lm}\}$ and then 
$\Psi_{\io}(B_1(\lm))=\{(\cd, 0,0)\}$. 
If the length of $w$ is equal to 1, then 
we can set $w=s_{i_1}$. As we have mentioned above,
the image of $B_{s_{i_1}}(\lm)$ by $\Psi_{\io}^{(\lm)}$ 
is 
$$
\{(\cd,x_2,x_1)\in {\mbox Im}(\Psi^{(\lm)}_{\io})\,|\,
x_2=x_3=\cd=0\}=\Sigma_{s_{i_1}}[\lm].
$$
Fix $w=s_{i_L}s_{i_{L-1}}\cd s_{i_2}s_{i_1}$ (reduced expression),
and set $w':=s_{i_{L-1}}\cd s_{i_2}s_{i_1}$.
By the hypothesis of the induction, we have
\begin{equation} 
\label{L-1}
\Psi_{\io}^{(\lm)}(B_{w'}(\lm))
=\Sigma_{w'}[\lm]:=\{(\cd,x_k,\cd,x_2,x_1)\in {\mbox Im}(\Psi^{(\lm)}_{\io})\,|\,
x_k=0{\mbox { for }}k>L-1\}.
\end{equation}
Here we show 
\begin{equation}
\label{con}
\til f_{i_L}^lx\in \Sigma_{w}[\lm]\cup\{0\}
\q({\hbox{for any $x\in \Sigma_{w'}[\lm]$ and any $l\in \ZZ_{\geq 0}$}}),
\end{equation}
by the induction on $l$.
For $x\in \Sigma_{w'}[\lm]$ and $k>L$ such that $i_k=i_L$, 
we have $\sigma_k(x)=\sigma_L(x)=0$ 
(as for $\sigma_k$ see (\ref{sigma})).
It follows from  (\ref{action-f}),
that if  $\til f_{i_L}x\ne0$, then 
its $k$-th entry is equal to 0.
Thus, we have
$$
\til f_{i_L}x\in \Sigma_{w}[\lm]\cup\{0\}.
$$
Suppose that 
\begin{equation}
\til f_{i_L}^lx\in \Sigma_{w}[\lm]
\label{step1}
\end{equation}
 and set its $L$-th entry $x'_L(\geq0)$.
By (\ref{step1}), we have $\sigma_k( \til f_{i_L}^lx)=0$
($k>L$ and $i_k=i_L$) 
and also we have $\sigma_L(\til f_{i_L}^lx)=x'_L\geq 0$.
This implies 
\begin{equation} 
\sigma_k(\til f_{i_L}^lx)\leq \sigma_L(\til f_{i_L}^lx)
\leq \sigma^{(i_L)}(\til f_{i_L}^lx).
\label{k-L}
\end{equation}
It follows from (\ref{action-f}) again
that we have 
$\til f_{i_L}^{l+1}x\in \Sigma_w[\lm]\cup\{0\}$ and then
$\Psi^{(\lm)}_{\io}(B_w(\lm))\subset \Sigma_w[\lm]$.

Next, we are going to show the opposite inclusion.
For any $x=(\cd,x_k,\cd,x_2,x_1)\in \Sigma_w[\lm]$, 
by (\ref{wt-vep-vp})
we have
\begin{equation}
\vep_{i_L}(x)={\mbox max}_{k;i_k=i_L}
\{\sigma_k(x),\sigma^{(i_L)}_0(x)\}
\geq \sigma_L(x)=x_L.
\label{x_L}
\end{equation}
Since the action of $\eit$ only reduces some entry in $x$,
we have 
$\til e_{i_L}^{\vep_{i_L}(x)}(x)\in \Sigma_w[\lm],$
(note that $\til e_{i_L}^{\vep_{i_L}(x)}(x)$ is never 0)
and 
\begin{equation}
\vep_{i_L}(\til e_{i_L}^{\vep_{i_L}(x)}(x))=0.
\label{ep-L}
\end{equation}
By (\ref{x_L}) and (\ref{ep-L}), we have
\begin{equation}
(\til e_{i_L}^{\vep_{i_L}(x)}(x))_L(=L{\hbox {-th entry of }}
\til e_{i_L}^{\vep_{i_L}(x)}(x))=0.
\label{L-zero}
\end{equation}
Thus, we have
$\til e_{i_L}^{\vep_{i_L}(x)}(x)\in \Sigma_{w'}[\lm]$.
Therefore, by Theorem \ref{kas}(ii), we get
\begin{equation}
x \in \til f_{i_L}^{\vep_{i_L}(x)}\Sigma_{w'}[\lm]
 =  \til f_{i_L}^{\vep_{i_L}}\Psi^{(\lm)}_{\io}(B_{w'}(\lm))
 \subset  \Psi^{(\lm)}_{\io}(B_w(\lm)).
\end{equation}
Now we obtain the opposite inclusion
$\Sigma_w[\lm]\subset \Psi^{(\lm)}_{\io}(B_w(\lm))$
and then completed the proof.\qed

Practically, we need the assumption ``ample''.
If $(\lm,\io)$ is ample, we can write Proposition \ref{polytope} in
the following form:
\begin{pro3}
If $(\lm,\io)$ is ample, we have
\begin{equation}
\Psi^{(\lm)}_{\io}(B_w(\lm))
(=\Sigma_w[\lm])=\{x=(x_k)\in \ZZ^{\ify}_{\io}[\lm]\,|\,
\begin{array}{l}
\vp(x)\geq 0
{\mbox{ for any }}\vp\in \Xi_{\io}[\lm],\,\\
x_k=0{\mbox{ for }}k>L.
\end{array}
\},
\label{ample-poly}
\end{equation}
where $\Xi_{\io}[\lm]$ is given in (\ref{Xi}).
\end{pro3}
Now, we obtain the convex ``polytope'' for $B_w(\lm)$.

\vskip5pt
In \cite{K3}, Kashiwara also introduced the crystal 
$B_w(\ify)\subset B(\ify)$. 
This possesses the following remarkable property:

\vskip5pt
\nd
If $b\in B(\ify)$ and $w\in W$ satisfy $\fit b\in B_w(\ify)$,
then $\fit ^kb\in B_w(\ify)$ for any $k\geq0$.

\vskip5pt
\nd 
This is used for proving Theorem \ref{kas} (iii).

It is characterized by the following;
\begin{enumerate}
\item
$B_w(\ify)=\{u_{\ify}\}$ if $w=1$.
\item
if $s_iw<w$, then
$B_w(\ify)=\bigcup_{k\geq 0}\fit^kB_{s_iw}(\ify).$
\end{enumerate}
This implies that 
$B_w(\ify)$ has also the similar description to $B_w(\lm)$.
\begin{pro3}
\begin{enumerate}
\item We have
$$
\Psi_{\io}(B_w(\ify))=\{(\cd,x_k,\cd,x_2,x_1)\in
{\rm Im}(\Psi_{\io})\,|\,
x_k=0{\mbox{ for }} k>L\}
$$
\item
If $\io$ satisfies the condition (\ref{posi}), we have
$$
\Psi_{\io}(B_w(\ify))=\left\{x=(\cd,x_k,\cd,x_2,x_1)\in\ZZ^{\ify}_{\io}\,
\left|\right.\,
\begin{array}{l}
\vp(x)\geq0{\mbox{ for any }}\vp\in \Xi_{\io}\\
x_k=0{\mbox{ for }} k>L
\end{array}
\right\},
$$
where $\Xi_{\io}$ is given in (\ref{Xi_io}).
\end{enumerate}
\end{pro3}

\subsection{Semi-simple cases}

In this subsection, let $\ge$ be a semi-simple Lie algebra,
$W$ be the corresponding Weyl group and $w_0\in W$ be the 
longest element with the length $l_0$.

In \cite[Proposition 4.2]{N3} 
we have shown by using 
the braid-type isomorphisms
that $B(\lm)$ can be embedded in the 
finite rank $\ZZ$-lattice $\ZZ^{l_0}$.
Here we obtain its simpler proof as an application of 
Proposition \ref{polytope}.
Indeed, in this case, since $V(\lm)=V_{w_0}(\lm)$,
we have $B(\lm)=B_{w_0}(\lm)$.
This implies:
\begin{pro3}
There exists the following embedding,
\begin{equation}
\Psi^{(\lm)}_{\io}:B(\lm)(=B_{w_0}(\lm))\mapright{\sim}
\Sigma_{w_0}[\lm]\hookrightarrow \ZZ^{l_0},
\end{equation}
where $\io$ is an infinite sequence of indices such that 
its first $l_0$ subsequence 
$i_{l_0},i_{l_0-1},\cd,i_1$ is a reduced longest word associated with 
the longest element $w_0$ $($see \cite[4.2]{N3}$)$.
\end{pro3}

\section{Extremal vectors}
\setcounter{equation}{0}
\renewcommand{\thesection}{\arabic{section}}
\renewcommand{\theequation}{\thesection.\arabic{equation}}

We still keep the notations of 3.2 and 
we do not necessarily assume that 
$(\lm,\io)$ is ample.

\subsection{Explicit description of extremal vectors}

For $w\in W$, we call $w\lm$ the extremal weight of 
$B(\lm)$ and call the unique element $u_{w\lm}\in B(\lm)_{w\lm}$ 
{\it extremal vector} with the extremal weight $w\lm$.

The image of $u_{w\lm}$ by $\Psi^{(\lm)}_{\io}$ 
is included in $\Sigma_w[\lm]$. We are going to 
determine it by the following way.

\newtheorem{pro4}{Proposition}[section]
\begin{pro4}
\label{extremal}
For $w\in W$ (length$(w)=L$),
set $x_w=(\cd,x_k,\cd,x_2,x_1):=\Psi^{(\lm)}_{\io}(u_{w\lm})$.
Then the element $x_w$ is given as the unique solution
of the following system of linear equations:
\begin{equation}
\left\{
\begin{array}{ll}
x_k=0 & {\mbox{ for }}k>L,\\
\beta^{(-)}_k(x)=0 & {\mbox{ for }}k\leq L,
\end{array}
\right.
\end{equation}
where the linear function $\beta^{(-)}_k$ is as in (\ref{beta--}).
\end{pro4}

\vskip 5pt 
{\sl Proof.\,}
The equations eq($L$)
$$
\beta^{(-)}_1(x)=\beta^{(-)}_2(x)=\cd
=\beta^{(-)}_L(x)=0.
$$
is the system of the linear equations in 
indeterminates $x_1,x_2,\cd,x_L$.
If we write eq($L$) in a matrix form 
$A\vec x=\vec\xi$ where $\vec x={}^t(x_1,\cd,x_L)$, 
due to the explicit form of $\beta^{(-)}_k$ in (\ref{beta--}),
the matrix $A$ 
is a triangular integer matrix whose diagonal entries are all 1 and 
the vector $\vec\xi={}^t(\xi_1,\cd,\xi_L)$ is given by
$\xi_k=\lan h_{i_k},\lm\ran$ if $\km=0$ and 
 otherwise $\xi_k=0$.
Thus, the equation eq($L$) can be solved uniquely
and all the entries are integers. We set the solution $(y_1,\cd,y_L)$.
Therefore, it suffices to show
$x_w(:=\Psi^{(\lm)}_{\io}(u_{w\lm}))=(\cd,0,0,y_L,\cd,y_1)$.
Let us show this by the induction on the length of $w$.
Set 
$w:=s_{i_L}s_{i_{L-1}}\cd s_{i_2}s_{i_1}$,
 $w':=s_{i_{L-1}}\cd s_{i_2}s_{i_1}$,
$y_w:=(\cd,0,0, y_L,y_{L-1},\cd,y_2,y_1)$ and 
$y_{w'}:=(\cd,0,0,y_{L-1},\cd,y_2,y_1)$.
Note that $y_{w'}$ is the unique solution of eq($L-1$)
and also the image of $x_{w'}$ by $\Psi^{(\lm)}_{\io}$
from the hypothesis of the induction.
Here we show the following lemma:
\newtheorem{lem4}[pro4]{Lemma}
\begin{lem4}
\label{ext-ext}
For $w$ and $w'$ as above, let $u_{w\lm}$ and $u_{w'\lm}$ be 
the corresponding extremal vectors. Then we have
\begin{equation}
u_{w\lm}=\til f_{i_L}^{\rm{max}} u_{w'\lm},
\label{ext-siki}
\end{equation}
where $\fit^{\rm max}u:=\fit^{\vp_i(u)}u$.
\end{lem4}

{\sl Proof of Lemma \ref{ext-ext}.\,}
By the definition of $\fit^{\rm max}$, 
$\til f_{i_L}^{\rm{max}} u_{w'\lm}\ne 0$.
Owing to the uniqueness of the extremal vector,
it suffices to show 
\begin{equation}
wt(\til f_{i_L}^{\rm{max}} u_{w'\lm})=w\lm.
\label{wt-w'}
\end{equation}
Let $S_L$ be the $i_L$-string in $B(\lm)$ including $u_{w'\lm}$.
By Theorem \ref{kas} (iii), we know that
$S_L\cap B_{w'}(\lm)$ is equal to 
 (1) $S_L$ or (2) $\{$highest weight vector in
$S_L$$\}$.
In the case (1), $u_{w'\lm}$ is the lowest weight vector in $S_L$
since $\mu\succ w'\lm$ 
when $B_{w'}(\lm)_{\mu}\ne\emptyset$.
This implies $\lan h_{i_L},w'\lm\ran\leq 0$.
Suppose that $\lan h_{i_L},w'\lm\ran<0$. 
Then we have 
\begin{equation}
w\lm=s_{i_L}(w'\lm)=w'\lm-\lan h_{i_L},w'\lm\ran\al_{i_L}\succ w'\lm,\nn
\end{equation}
which contradicts Theorem \ref{kas} (ii).
Thus, in this case we have $\lan h_{i_L}w'\lm\ran=0$ 
and then  the length of $S_L=0$.
This means $\vp_{i_L}(u_{w'\lm})=0$ and then $u_{w\lm}=u_{w'\lm}$.
In the case (2), since $\vep_{i_L}(u_{w'\lm})=0$, we have 
$\vp_{i_L}(u_{w'\lm})=\lan h_{i_L},wt(u_{w'\lm})\ran$ and then
\begin{eqnarray*}
wt(\til f_{i_L}^{\rm{max}} u_{w'\lm})
&=& w'\lm-\vp_{i_L}(u_{w'\lm})\al_{i_L}\\
&=& w'\lm-\lan h_{i_L},wt(u_{w'\lm})\ran\al_{i_L}\\
&=& s_{i_L}(w'\lm)=w\lm.
\end{eqnarray*}
Now, we obtain (\ref{wt-w'}) and then completed the proof of Lemma
\ref{ext-ext}
\qed

By this lemma, it suffices to show
\begin{equation} 
\til f_{i_L}^{\rm{max}} y_{w'}=y_w,
\label{yy}
\end{equation}
Let us see how $\til f_{i_L}$ acts on $y_{w'}$.
For $k\in\ZZ_{\geq1}$ and $m\in\ZZ_{>0}$ we set 
$k^{(\pm1 )}:=k^{(\pm)}$, 
$k^{(+m)}:=(k^{(+(m-1))})^{(+)}$ and 
$k^{(-m)}:=(k^{(-(m-1))})^{(-)}$ 
(as for $k^{(\pm)}$, see 2.3).
For $m\geq 1$, we have
\begin{equation}
\sigma_{L^{(-m)}}+\beta^{(-)}_{L^{(-m)}}=\sigma_{L^{(-m-1)}}.
\label{sig-sig}
\end{equation}
Since $\beta^{(-)}_k(y_{w'})=0$ for $k\leq L-1$, we have by (\ref{sig-sig})
\begin{equation}
\sigma_{L^{(-)}}(y_{w'})=
\sigma_{L^{(-2)}}(y_{w'})=\cd
=\sigma_{L^{(-m)}}(y_{w'})=\cd
=\sigma_0^{(i_L)}(y_{w'}).
\label{sigs}
\end{equation}
Now we consider the following cases:

(a) $\beta^{(-)}_L(y_{w'})<0$.\qq (b) $\beta^{(-)}_L(y_{w'})\geq0$.

\nd
In the case (a), we have
\begin{equation}
\begin{array}{lll}
0=&\cd =& \sigma_{L^{(+m)}}(y_{w'})=\cd = \sigma_L(y_{w'})\\
&&>
\sigma_{L^{(-)}}(y_{w'})=
\cd=\sigma_{L^{(-m)}}(y_{w'})=\cd
=\sigma_0^{(i_L)}(y_{w'}).
\label{case1}
\end{array}
\end{equation}
It follows from (\ref{action-f}) that 
\begin{equation}
\til f_{i_L}(y_{w'})=(\cd,1,y_{L-1},\cd,y_1).
\label{f-il}
\end{equation}
In the case (b), we have
\begin{equation}
\begin{array}{lll}
0=&\cd &=\sigma_{L^{(+m)}}(y_{w'})=\cd = \sigma_L(y_{w'})\\
&\leq&
\sigma_{L^{(-)}}(y_{w'})=
\cd=\sigma_{L^{(-m)}}(y_{w'})=\cd
=\sigma_0^{(i_L)}(y_{w'}),
\end{array}
\label{*}
\end{equation}
which implies $\til f_{i_L}(y_{w'})=0$ by (\ref{action-f}).
In this case, by (\ref{wt-vep-vp-1}) and (\ref{wt-vep-vp})  we have
\begin{equation}
\vep_{i_L}(y_{w'})=\sigma_0^{(i_L)}(y_{w'})=-\lan
h_{i_L},wt(y_{w'})\ran,
\label{wt-yw'}
\end{equation}
and then also by (\ref{wt-vep-vp-3}) we have
\begin{equation}
\vp_{i_L}(y_{w'})=\vep_{i_L}(y_{w'})+\lan h_{i_L},wt(y_{w'})\ran=0,
\end{equation}
which implies that $y_{w'}$ is the lowest weight vector 
in the $i_L$-string including $y_{w'}$. 
Thus, the case (b) corrsponds to the case (1) 
in the proof of Lemma \ref{ext-siki}.
So the length of $i_L$-string is 0 and then 
we have 
$\vep_{i_L}(y_{w'})=0$.
By (\ref{*}) and (\ref{wt-yw'}) we have
$0=\sigma^{(i_L)}_0(y_{w'})=\cd =\sigma_{L^{(-)}}(y_{w'})=
\sigma_L(y_{w'})=0$, and then
$\beta^{(-)}_L(y_{w'})0$ by (\ref{sig-sig}),
which means $y_w=y_{w'}(=\til f_{i_L}^{\rm max}y_{w'})$, that is,
(\ref{yy}) with $y_L=0$.

In the case (a), we can suppose that $\vp_{i_L}(y_{w'})>0$. 
Let us show 
\begin{equation}
\til f_{i_L}^k(y_{w'})=(\cd,0,0,k,y_{L-1},\cd,y_2,y_1),
\label{f-k}
\end{equation}
for $1\leq k\leq \vp_{i_L}(y_{w'})$ by the induction on $k$.
Assuming (\ref{f-k}) ($1\leq k<\vp_{i_L}(y_{w'})$), let us see
$\til f_{i_L}^{k+1}(y_{w'})$.
Set $\bar y:=\til f_{i_L}^k(y_{w'})
=(\cd,0,0,k,y_{L-1},\cd,y_2,y_1)$.
In the case (a), by the argument for the case (1) 
in the proof of Lemma \ref{ext-ext},
 we have
$\vep_{i_L}(y_{w'})=0$ and then by 
$\sigma_L(y_{w'})=0$ and (\ref{**}), 
\begin{equation}
\vp_{i_L}(y_{w'})=\lan h_{i_L},wt(y_{w'}) \ran
=-\sigma^{(i_L)}_0(y_{w'})=\cd =-\sigma_{L^{(-)}}(y_{w'})
=-\beta^{(-)}_{L}(y_{w'})
\label{y-1}
\end{equation}
Set $F:=\vp_{i_L}(y_{w'})=-\beta^{(-)}_{L}(y_{w'})$.
On the other hand, if $i_l=i_L$, 
we have
\begin{equation}
\sigma_l(\bar y)=\sigma_l(y_{w'})+2k =-F+2k\q
(l<L),\qq
\sigma_{L}(\bar y)=k.
\label{y-2}
\end{equation}
It follows from (\ref{y-1}) and (\ref{y-2}) that
if $k<F$, we have $\sigma_{L^{(-n)}}(\bar y)=2k-F<k
=\sigma_L(\bar y)\geq0=\sigma_{L^{(+m)}}(\bar y)$ 
$(m,n>0)$ and then
by (\ref{action-f})
\begin{equation}
\til f_{i_L}(\bar y)=(\cd,0,0,k+1,y_{L-1},\cd,y_1).
\end{equation}
Hence, we obtain
\begin{equation}
\til f^{\rm{max}}_{i_L}(y_{w'})=
\til f_{i_L}^{F}(y_{w'})=
(\cd,0,0,-\beta^{(-)}_{L}(y_{w'}),y_{L-1},\cd,y_1).
\end{equation}
Here $(y_1,\cd,y_{L-1},-\beta^{(-)}_L(y_{w'}))$
satisfies the equations eq($L$) and then it follows that 
$y_w=\til f^{\rm max}_{i_L}(y_{w'})$.
\qed

{\sl Remark.\,}
Note that we obtain Proposition \ref{extremal} without the 
assumption ``ample''.
In \cite[Example 3.9]{N3}, we introduced the 
``non-ample'' exmple: $\ge=A_3$ and $\io=212321$.
But, even in this case, applying Propsition \ref{extremal}
we have
\begin{eqnarray*}
&&x_{s_3s_2s_1} =(0,0,0,\lm_1+\lm_2+\lm_3,\lm_1+\lm_2,\lm_1),\\
&&x_{s_2s_3s_2s_1}=(0,0,\lm_3,\lm_1+\lm_2+\lm_3,\lm_1+\lm_2,\lm_1),\\
&&x_{s_1s_2s_3s_2s_1}=(0,\lm_2+\lm_3,
\lm_3,\lm_1+\lm_2+\lm_3,\lm_1+\lm_2,\lm_1),\\
&&x_{s_2s_1s_2s_3s_2s_1}=(\lm_2,\lm_2+\lm_3,
\lm_3,\lm_1+\lm_2+\lm_3,\lm_1+\lm_2,\lm_1),
\end{eqnarray*}
where $\lm_i=\lan h_i,\lm\ran$.

\subsection{Rank 2 cases}

We apply Propositon \ref{extremal} to arbitrary rank 2 cases.

First we review the result in \cite{N2}. 
The setting here is same as those in \cite{N2}.
We set $I=\{1,2\}$,
and $\io = (\cd,2,1,2,1)$.
The Cartan matrix is given by:
$$
\lan h_1,\al_1\ran= \lan h_2,\al_2\ran=2, \,\, \lan h_1,\al_2\ran=-c_1,
\,\, \lan h_2,\al_1\ran=-c_2.
$$
Here we either have $c_1 = c_2 = 0$, or both $c_1$ and $c_2$ are
positive integers.
We set $X = c_1 c_2 - 2$, and define the integer sequence
$a_l = a_l (c_1, c_2)$ for $l \geq 0$ by setting $a_0 = 0, \, a_1 = 1$
and, for $k \geq 1$,
\begin{equation}
a_{2k}  = c_1 P_{k-1} (X), \,\,
a_{2k+1} = P_k (X) + P_{k-1} (X),
\label{defcoeff}
\end{equation}
where the $P_k (X)$ are {\it Chebyshev polynomials} given 
by the following generating function:
\begin{equation}
 \sum_{k \geq 0} P_k (X) z^k = (1 - X z + z^2)^{-1}.
\label{gen-cheb}
\end{equation}
Here define $a'_l(c_1,c_2):=a_l(c_2,c_1)$.
The several first Chebyshev polynomials and terms $a_l$ are given by
$$P_0 (X) = 1, \, P_1 (X) = X, \,
P_2 (X) = X^2 - 1, \, 
P_3 (X) = X^3 - 2 X,$$
$$a_2 = c_1, \, a_3 = c_1 c_2 - 1, \, a_4 = c_1 (c_1 c_2 - 2),$$
$$a_5 = (c_1 c_2 - 1)(c_1 c_2 - 2) - 1, \, 
a_6 = c_1 (c_1 c_2 - 1)(c_1 c_2 - 3).$$
Let $l_{\rm max} = l_{\rm max} (c_1, c_2)$ be the minimal index
$l$ such that $a_{l+1} < 0$
(if $a_l \geq 0$ for all $l \geq 0$,
then we set $l_{\rm max} = + \infty$). 
By inspection, if $c_1 c_2 = 0$ (resp. $1,2,3$) 
then $l_{\rm max} = 2$ (resp. $3, 4, 6$).
Furthermore, if $c_1 c_2 \leq 3$ then $a_{l_{\rm max}} = 0$ and 
$a_l > 0$ for $1 \leq l < l_{\rm max}$.
On the other hand, if $c_1 c_2 \geq 4$, i.e., $X \geq 2$,
it is easy to see from (\ref{gen-cheb})
that $P_k (X) > 0$ for $k \geq 0$, hence 
$a_l > 0$ for $l \geq 1$; in particular, in this case
$l_{\rm max} = + \infty$.  

\begin{pro4}[\cite{N2}]
\label{rank 2-thm}
In the rank 2 case, for a dominant integral weight 
$\lm=m_1\Lm_1+m_2\Lm_2$ $(m_1,m_2\in \ZZ_{\geq0})$ 
the image of the embedding $\Psi^{(\lm)}_{\io}$ 
is given by
\begin{equation}
{\rm Im} \,(\Psi^{(\lm)}_{\io}) = \left\{(\cd,x_2,x_1)\in\ZZ_{\geq0}^{\ify}: 
\begin{array}{l}
x_k = 0 \,\,{\rm for}\,\,k > l_{\rm max},\,\,m_1\geq x_1, \\
a_l x_l -a_{l-1} x_{l+1} \geq0 ,\\
m_2+a'_{l+1}x_l-a'_lx_{l+1}\geq0,\\
{\rm for}\,\,1 \leq l < l_{\rm max}
\end{array}
\right\}.
\label{rank 2-poly}
\end{equation}
\end{pro4}

Note that the cases when $l_{\rm max} < +\infty$, or equivalently, 
the image ${\rm Im} \,(\Psi_{\io}^{(\lm)})$ 
is contained in a lattice of finite rank,
just correspond to the Lie algebras 
$\ge=$ $A_1 \times A_1$, $A_2$, $B_2$ or $C_2$, $G_2$.

For $L\in \ZZ_{\geq0}$ $(L\leq l_{\rm max})$,  we set 
\begin{equation}
w_L:=\left\{
\begin{array}{ll}
s_1(s_2s_1)^l & {\mbox{ if }} L=2l+1,\\
(s_2s_1)^l    & {\mbox{ if }} L=2l.\\
\end{array}
\right.
\end{equation}
For a dominant integrable weight $\lm=m_1\Lm_1+m_2\Lm_2\in P_+$,
we define the integer sequence $d_k=d_k(c_1,c_2)$ ($k\geq 1$) as follows:
\begin{equation}
d_k(c_1,c_2):=m_1a_k(c_2,c_1)+m_2a_{k-1}(c_1,c_2),
\end{equation}
where $\{a_k\}_{k\geq0}$ is the integer sequence given in (\ref{defcoeff}).
\begin{pro4}
The image $x_{w_L}$ of the extremal vector $u_{w_L\lm}\in B_{w_L}(\lm)$ 
$(\lm=m_1\Lm_1+m_2\Lm_2\in P_+)$ associated 
with $w_L\in W$ can be described as follows:
\begin{equation}
x_{w_L}\left(:=\Psi^{(\lm)}_{\io}(u_{w_L\lm})\right)
=(\cd,0,0,d_L,d_{L-1},\cd,d_2,d_1).
\label{r2-ex}
\end{equation}
\end{pro4}

{\sl Proof.\,}
In this setting, we have
\begin{eqnarray*}
& \beta^{(-)}_1(x)=x_1-m_1,\,\,
\beta^{(-)}_2(x)=x_2-c_2x_1-m_2,\\
&\beta^{(-)}_{2k+1}(x)=x_{2k+1}-c_1x_{2k}+x_{2k-1},\,\,
\beta^{(-)}_{2k+2}(x)=x_{2k+2}-c_2x_{2k+1}+x_{2k},\,\,
(k\geq1).
\end{eqnarray*}
We shall show $d_1,d_2,\cd,d_L$ are the solutions of the equations
\begin{equation}
\beta^{(-)}_1(x)=\beta^{(-)}_2(x)=\cd =\beta^{(-)}_L(x)=0.
\label{beta-beta}
\end{equation}
Solving $\beta^{(-)}_1(x)=\beta^{(-)}_2(x)=0$, we have
\begin{equation}
x_1=m_1=m_1a'_1+m_2a_0=d_1,\qq x_2=m_1c_2+m_2=m_1a'_2+m_2a_1=d_2.
\label{x_1x_2}
\end{equation} 
Here note that we can write $d_k=m_1a'_k+m_2a_{k-1}$.
By the definition of $a_k$, we can easily see that 
$\{a_k\}$ (resp. $\{a'_k\}$) is uniquely determined by 
\begin{eqnarray*}
&\hspace{-7pt} a_0=0,\,\,a_1=1,\,\,
a_{2k+1}=c_2a_{2k}-a_{2k-1},\,\,
a_{2k+2}=c_1a_{2k+1}-a_{2k},\,\,(k\geq0).\\
&({\rm resp.}\,\,
 a'_0=0,\,\,a'_1=1,\,\,
a'_{2k+1}=c_1a'_{2k}-a'_{2k-1},\,\,
a'_{2k+2}=c_2a'_{2k+1}-a'_{2k}.).\nn
\end{eqnarray*} 
Here for $k\geq1$ we have
$d_{2k+1}-c_1d_{2k}+d_{2k-1}
=m_1(a'_{2k+1}-c_1a'_{2k}+a'_{2k-1})+
m_2(a_{2k}-c_1a_{2k-1}+a_{2k-2})=0$
and
$d_{2k+2}-c_2d_{2k+1}+d_{2k}
=m_1(a'_{2k+2}-c_2a'_{2k+1}+a'_{2k})+
m_2(a_{2k+1}-c_2a_{2k}+a_{2k-1})=0$,
which implies that 
$(d_1,d_2,\cd,d_L)$ is the unique solution of (\ref{beta-beta}).
Now, we obtain the desired result.\qed

In conclusion of this section, 
we illustlate the case of $A^{(1)}_1$, that is, 
 $c_1 = c_2 = 2$.
In this case, $X = c_1 c_2 - 2 = 2$.
It follows at once from (\ref{gen-cheb}) that $P_k (2) = k+1$;
hence, (\ref{defcoeff}) gives $a_l = l$ for $l \geq 0$.
We see that for type $A^{(1)}_1$, 
$$
x_{w_L}=(\cd,0,0,Lm_1+(L-1)m_2,\cd,km_1+(k-1)m_2,\cd,2m_1+m_2,m_1).
$$

\vspace{5pt}

{\bf Acknoledgements\,}
The author would like to acknowledge M.Okado
for his interest to my work.
Indeed, this work was motivated by his question
about the relation of polyhedral realizations and 
the crystals of Demazure modules.

\end{document}